\documentclass{tran-l}

\def\tratto{\mbox{\rule{2mm}{.2mm}$\;\!$}}
\def\split{\genfrac{}{}{0pt}1}

\theoremstyle{plain}
\newtheorem{Theorem}{Theorem}[section]
\newtheorem{Lemma}[Theorem]{Lemma}
\newtheorem{Corollary}[Theorem]{Corollary}
\newtheorem{Proposition}[Theorem]{Proposition}
\newtheorem{Conjecture}[Theorem]{Conjecture}

\theoremstyle{definition}

\newtheorem{Example}[Theorem]{Example}

\theoremstyle{remark}
\newtheorem{Remark}[Theorem]{Remark}

\begin{document}

\title{Core and residual intersections of ideals}

\author[A. Corso]{Alberto Corso}
\address{Department of Mathematics, University of Kentucky, Lexington,
Kentucky 40506} \email{corso@ms.uky.edu}
\urladdr{www.ms.uky.edu/{\textasciitilde}corso}

\author[C. Polini]{Claudia Polini}
\address{Department of Mathematics, University of Oregon,
Eugene, Oregon 97403} \curraddr{Department of Mathematics,
University of Notre Dame, Notre Dame, Indiana 46556}
\email{cpolini@nd.edu}
\urladdr{www.nd.edu/{\textasciitilde}cpolini}

\author[B. Ulrich]{Bernd Ulrich}
\address{Department of Mathematics, Michigan State University,
East Lansing, Michigan 48824} \curraddr{Department of Mathematics,
Purdue University, West Lafayette, Indiana 47907}
\email{ulrich@math.purdue.edu}
\urladdr{www.math.purdue.edu/{\textasciitilde}ulrich}

\dedicatory{To Professor Craig Huneke on the occasion of his
fiftieth birthday}

\date{April 10, 2001}

\thanks{The first author was partially supported by the NATO/CNR Advanced
Fellowships Programme during an earlier stage of this work. The
second and third authors were partially supported by the NSF}



\begin{abstract}
D. Rees and J. Sally defined the core of an $R$-ideal $I$ as the
intersection of all $($minimal$)$ reductions of $I$. However, it
is not easy to give an explicit characterization of it in terms of
data attached to the ideal. Until recently, the only case in which
a closed formula was known is the one of integrally closed ideals
in a two-dimensional regular local ring, due to C. Huneke and I.
Swanson. The main result of this paper explicitly describes the
core of a broad class of ideals with good residual properties in
an arbitrary local Cohen--Macaulay ring. We also find sharp bounds
on the number of minimal reductions that one needs to intersect to
get the core.
\end{abstract}

\maketitle

\section{Introduction}

Let $R$ be a Noetherian local ring and $I$ one of its ideals. D.
Rees and J. Sally defined in \cite{RS} the {\it core} of $I$,
denoted ${\rm core}(I)$, to be the intersection of all
$($minimal$)$ reductions of $I$. Roughly, a reduction can be
thought of as a simplification of the original ideal: It is a
notion introduced by D.G. Northcott and D. Rees \cite{NR} and has
recently played a crucial role in the study of Rees algebras of
ideals. Hence, ${\rm core}(I)$ can be viewed as the
`simplification of all simplifications' of $I$! The core of $I$
naturally appears also in the context of the celebrated
Brian\c{c}on--Skoda theorem, which -- in one of its easiest
formulations -- says that if $R$ is a regular local ring of
dimension $d$ and $I$ is an ideal then $\overline{I^d} \subset
{\rm core}(I)$, where $\overline{I^d}$ denotes the integral
closure of $I^d$. Our goal can be framed in terms of giving an
explicit description of ${\rm core}(I)$ from data attached to the
ideal $I$. Without further assumptions, our goal is ambitious --
if not hopeless -- as ${\rm core}(I)$ is the intersection of an
{\it a priori} infinite number of ideals. Thus, our first task is
to identify a natural setting where our calculations will go
through: Residual intersections of ideals provide such a
framework.

To be more specific, an ideal $I$ is said to be {\it integral}
over an ideal $J \subset I$ if the inclusion of Rees algebras
${\mathcal R}(J) \hookrightarrow {\mathcal R}(I)$ is module
finite. The {\it integral closure} $\overline{I}$ of $I$ is then
defined to be the largest ideal integral over $I$ and the ideal
$I$ is {\it integrally closed} $($or {\it complete}$)$ if
$\overline{I}=I$. Alternatively, if $I$ is integral over $J$ the
ideal $J$ is  said to be a {\it reduction} of $I$. Equivalently,
$J$ is a reduction of $I$ if $I^{r+1}=JI^r$ for some non negative
integer $r$. The least such $r$ is called the reduction number
$r_J(I)$ of $I$ with respect to $J$. One then defines the {\it
reduction number} of $I$ to be the least $r_J(I)$, where $J$
varies among all minimal reductions of $I$. A reduction is said to
be {\it minimal} if it is minimal with respect to containment.
Minimal reductions of an ideal have the pleasing property of
carrying most of the information about the original ideal but, in
general, with fewer generators. If the residue field $k$ of the
ring is infinite the minimal number of generators of a minimal
reduction equals the {\it analytic spread} $\ell$ of the ideal
$I$, i.e., the dimension of the {\it special fiber ring}
${\mathcal R}(I) \otimes_R k$ of $I$.

Little is known about the structure and the properties of ${\rm
core}(I)$. In \cite[2.6]{RS} Rees and Sally proved that if $I$ is
an ${\mathfrak m}$-primary ideal of a $d$-dimensional local
Cohen--Macaulay ring $(R, {\mathfrak m})$ and ${\mathcal A}$ is a
`universal $d$-generated ideal' in $I$ then ${\mathcal A} \cap R
\subset {\rm core}(I)$. Later, \cite[3.9]{HS} provided the first
work in the literature with a `closed formula' for the core of an
ideal. More precisely, C. Huneke and I. Swanson showed that the
core of integrally closed ideals in two-dimensional regular local
rings is still integrally closed and is given by a formula that
involves an ideal of minors of any presentation matrix of the
ideal. They also showed that the core of such ideals is closely
related to the adjoint, introduced by J. Lipman in \cite{L}, and
proved several arithmetical properties of the core. The techniques
used involve a systematic application of the Hilbert--Burch
theorem, the notion of quadratic transforms and Rees valuations.
More recently, we showed in \cite{CPU}, under fairly general
assumptions, for which we refer the reader to \cite{CPU}, that:
\begin{itemize}
\item
${\rm core}(I)$ is the intersection of finitely many general
reductions of $I$;
\item
${\rm core}(I)$ is the contraction to $R$ of a `universal ideal in
$I$';

\item
${\rm core}(I)$ behaves well under flat extensions.
\end{itemize}
From the second result we also deduced an expression for ${\rm
core}(I)$ as a colon ideal in a polynomial ring over $R$ that
allows -- at least in principle -- for explicit calculations.

The spirit of \cite{CPU} was close to the one of \cite{RS}. In
this paper we shift interest instead: Our main goal is to give an
explicit formula for ${\rm core}(I)$ in the spirit of \cite{HS},
by which we mean a formula that only involves operations inside
the ring $R$ itself. The ideals under consideration have slightly
more structure than the ones studied in \cite{CPU}. Nevertheless,
our results substantially extend the ones of Huneke and Swanson.
The techniques we use are rather different and require some of the
machinery developed in \cite{CPU}. To arrive at our more general
formula we needed to observe that residual intersections are the
correct objects to replace the Fitting ideals occurring in the
Huneke--Swanson result. Thus we can describe the core of ideals
$I$ that are {\it balanced} $($i.e., $J \colon I$ is independent
of the minimal reduction $J$ of $I$, see \cite{U2}$)$ or have the
{\it expected reduction number} $\leq \ell-g+1$, where $g = {\rm
ht}\, I$. Interestingly enough these conditions also turn out to
be necessary for the validity of our formula!

Section~2 contains our main results and their proofs. We let $(R,
{\mathfrak m})$ be a local Cohen--Macaulay ring and we focus on
ideals that satisfy specific bounds on the number of local
generators up to a certain codimension $($property $G_{\ell})$ and
have good residual $S_2$ properties. We first prove that for this
kind of ideals the property of being balanced is equivalent to the
inclusion $(J \colon I)I \subset {\rm core}(I)$, where $J$ is any
minimal reduction of $I$ $($see Proposition~\ref{thm1}$)$. We then
devote most of the section to showing in Theorem~\ref{converse}
that this inclusion is indeed an equality, that is
\[
{\rm core}(I) = (J \colon I) I,
\]
for any minimal reduction $J$ of $I$.  If in addition $R$ is
Gorenstein, $\ell \geq 1$, and ${\rm depth}\, R/I^j \geq {\rm
dim}\, R/I -j+1$ for $1 \leq j \leq \ell-g+1$, we actually prove
that the above formula is equivalent to the ideal $I$ having the
expected reduction number. Under these assumptions, we also obtain
that ${\rm core}(I)$ is the intersection of $\ \displaystyle \ell
\cdot {\mu(I)-g \choose \ell-g+1} +1 \ $ general minimal
reductions -- a bound that is sharp in several cases of interest.
To prove Theorem~\ref{converse} we first make a delicate use of
residual intersection techniques to reduce to the case of an
${\mathfrak m}$-primary ideal in a one-dimensional local
Cohen--Macaulay ring. We then assemble several intermediate
results and lemmas proved earlier as well as a key idea of Huneke
and Swanson. We finish the section with an application of
Theorem~\ref{converse} which says that the core of a normal
balanced ideal is integrally closed $($see
Theorem~\ref{normality}$)$. However, there are examples showing
that the integral closedness of $I$ alone is not sufficient to
guarantee the one of ${\rm core}(I)$ $($see
Example~\ref{bad-pfaff}$)$.

Section~3 lists various classes of ideals for which we can
explicitly describe the core from a matrix presenting the ideal.
They include the case of perfect ideals of height two or perfect
Gorenstein ideals of height three $($see
Corollaries~\ref{explicit1},~\ref{corHS},~\ref{explicit2}
and~\ref{explicit3}$)$. In particular, Theorem~\ref{converse}
recovers \cite[3.9]{HS} $($see Corollary~\ref{corHS}$)$.

In Section~4 we study the core of links of symbolic powers $($see
Proposition~\ref{symbolic}$)$ and the core of powers of complete
intersections $($see Proposition~\ref{proposition4.4}$)$. This
result applies to even symbolic powers of self-linked perfect
ideals of height $2$ and we are able to show that their core is an
odd symbolic power of the ideal $($see
Corollary~\ref{even-core}$)$. Self-linked ideals have interested
many researchers including D. Ferrand, M. Kumar, L. Szpiro, G.
Valla and more recently J. Herzog and B. Ulrich \cite{HeU}, and S.
Kleiman and B. Ulrich \cite{KU}. In particular it was shown in
\cite{KU} that these ideals are, roughly, in correspondence with
Gorenstein perfect algebras of grade one that are birational onto
their image -- a result that was inspired by recent work in the
theory of central projections onto hypersurfaces.

Finally, in Section 5 we discuss Conjecture~\ref{conj}: It asks
whether -- for any non-nilpotent ideal $I$ with $G_{\ell}$ and
good residual properties -- the core is given by the formula
\[
{\rm core}(I) = (J^r \colon I^r)I = (J^r \colon I^r)J = J^{r+1}
\colon I^r,
\]
for any minimal reduction $J$ of $I$ with reduction number $r$.
The point is to move away from the balancedness condition, which
was required in Theorem~\ref{converse}. This formula has solid
theoretical foundation and an extensive computer evidence in its
support. Moreover it is sharp, in the sense that \cite[4.11]{CPU}
provides a counterexample if the assumptions are relaxed. We show
in Proposition~\ref{agreement} that this formula agrees with the
one obtained in Theorem~\ref{converse}, at least if we add some
additional assumptions on the ring $R$ and the powers of the ideal
$I$.

\section{Residually $S_2$ ideals}

We begin by reviewing some facts and results from \cite{CEU} and
\cite{CPU}. Let $R$ be a Noetherian ring, $I$ an $R$-ideal of
height $g$, and $s$ an integer. Recall that $I$ satisfies
condition $G_s$ if for every prime ideal ${\mathfrak p}$
containing $I$ with ${\rm dim}\, R_{\mathfrak p} \leq s-1$, the
minimal number of generators $\mu(I_{\mathfrak p})$ of
$I_{\mathfrak p}$ is at most ${\rm dim}\, R_{\mathfrak p}$. The
ideal $I$ is said to have property $G_{\infty}$ if $G_s$ holds for
every $s$. A proper $R$-ideal $K$ is called an {\it $s$-residual
intersection} of $I$, if there exists an $s$-generated ideal
${\mathfrak a} \subset I$ so that $K= {\mathfrak a} \colon I$ and
the height of $K$ is at least $s \geq g$. If in addition the
height of $I + K$ is at least $s+1$, then $K$ is said to be a {\it
geometric $s$-residual intersection} of $I$. We say that $I$
satisfies $AN_s$ $(AN_s^{-})$ if $R/K$ is Cohen--Macaulay for
every $i$-residual intersection $($geometric $i$-residual
intersection, respectively$)$ $K$ of $I$ and every $i \leq s$. The
ideal $I$ is called {\it $s$-residually $S_2$ $($weakly
$s$-residually $S_2)$} if $R/K$ satisfies Serre's condition $S_2$
for every $i$-residual intersection $($geometric $i$-residual
intersection, respectively$)$ $K$ of $I$ and every $i \leq s$.
Finally, whenever $R$ is local, we say $I$ is {\it universally
$s$-residually $S_2$ $($universally weakly $s$-residually $S_2)$}
if $IS$ is $s$-residually $S_2$ $($weakly $s$-residually $S_2)$
for every ring $S=R(x_1, \ldots, x_n)$ with $x_1, \ldots, x_n$
variables over $R$.

If $(R, {\mathfrak m})$ is a local Cohen--Macaulay ring of
dimension $d$ and $I$ an $R$-ideal satisfying $G_s$, then $I$ is
universally $s$-residually $S_2$ in the following cases:
\begin{itemize}

\item[$($1$)$]
$R$ is Gorenstein, and the local cohomology modules $H_{\mathfrak
m}^{d-g-j}(R/I^j)$ vanish for $1 \leq j \leq s-g+1$, or
equivalently, ${\rm Ext}_R^{g+j}(R/I^j, R) = 0$ for $1 \leq j \leq
s-g+1$ $($see \cite[4.1 and 4.3]{CEU}$)$.

\item[$($2$)$]
$R$ is Gorenstein, and ${\rm depth}\, R/I^j \geq {\rm dim}\, R/I -
j + 1$ for $1 \leq j \leq s-g+1$ $($see \cite[2.9$(${\it
a}$)$]{U}$)$.

\item[$($3$)$]
$I$ has sliding depth, which means that the $i^{\rm th}$ Koszul
homology modules $H_i$ of a generating set $f_1, \ldots, f_n$ of
$I$ satisfy ${\rm depth}\, H_i \geq {\rm dim}\, R - n + i$ \, for
every $i$ $($see \cite[3.3]{HVV}$)$.
\end{itemize}

In fact condition $($2$)$ implies $($1$)$ and the property $AN_s$
\cite[2.9$(${\it a}$)$]{U}. Assumptions $($2$)$ and $($3$)$ in
turn are satisfied by strongly Cohen--Macaulay ideals, i.e.,
ideals whose Koszul homology modules are Cohen--Macaulay. This
condition always holds if $I$ is a Cohen--Macaulay almost complete
intersection or a Cohen--Macaulay deviation two ideal of a
Gorenstein ring \cite[p.{\,}259]{AH}. It is also satisfied for any
ideal in the linkage class of a complete intersection
\cite[1.11]{H1}: Standard examples include perfect ideals of
height two and perfect Gorenstein ideals of height three.

\medskip

The next proposition contains our first result: It relates the
shape of the core of $I$ to the property of $I$ being balanced.

\begin{Proposition}\label{thm1}
Let $R$ be a local Cohen--Macaulay ring with infinite residue
field. Let $I$ be an $R$-ideal of height $g$ and analytic spread
$\ell$. Suppose that $I$ satisfies property $G_{\ell}$ and is
weakly $(\ell-1)$-residually $S_2$. Then the following conditions
are equivalent:
\begin{itemize}
\item[$(${\it a}$)$]
$(J \colon I) I \subset {\rm core}(I)$ for some minimal reduction $J$ of $I$;

\item[$(${\it b}$)$]
$(J \colon I) I \subset {\rm core}(I)$ for every minimal reduction $J$ of $I$;

\item[$(${\it c}$)$]
$J \colon I$ does not depend on the minimal reduction $J$ of $I$.
\end{itemize}
If any of these conditions holds then $I^2 \subset {\rm core}(I)$.
\end{Proposition}
\begin{proof}
We may assume that $\ell \not= \mu(I)$. Let $H$ be any
minimal reduction of $I$.

If $(${\it a}$)$ holds then $(J \colon I)I \subset {\rm core}(I)
\subset H$, hence $J \colon I \subset H \colon I$. Since $I$ is
$G_{\ell}$ and weakly $(\ell-1)$-residually $S_2$, $J \colon I$
and $H \colon I$ are unmixed ideals of the same dimension
$($\cite[2.1$(\/\/{\it f}\/)$]{CPU}$)$. On the other hand
$e(R/J\colon I) = e(R/H \colon I)$ by \cite[4.2 and 4.4]{CPU}.
Thus $J \colon I = H \colon I$, which gives $(${\it c}$)$. If
$(${\it c}$)$ holds then $(J \colon I)I = (H \colon I)I \subset H$
for every minimal reduction $J$ of $I$, and we obtain $(${\it
b}$)$. Obviously $(${\it b}$)$ implies $(${\it a}$)$.

Finally, if any of these conditions holds then by $(${\it c}$)$,
$I \subset H \colon I$ since $I$ is the sum of its minimal
reductions. Hence $I^2 \subset {\rm core}(I)$.
\end{proof}

Our next goal is to strengthen Proposition~\ref{thm1} by showing
that with essentially the same assumptions one actually has that
${\rm core}(I) = (J \colon I) I$ for any minimal reduction $J$ of
$I$. We obtain this result in Theorem~\ref{converse}: It requires
the assemblage of several facts that we are going to prove next.
The proof of Lemma~\ref{lemma-HS} below has been inspired by the
one of \cite[3.8]{HS}. To simplify our notation if $x, y$ are
elements of $R$ and $I \subset R$ then by $x \colon I$ and $x
\colon y$ we mean $(x) \colon I$ and $(x) \colon (y)$,
respectively.

\medskip

\begin{Lemma}~\label{lemma-HS}
Let $(R, {\mathfrak m}, k)$ be a Noetherian local ring. Let $K$ be
an $R$-ideal and let $x$ and $y$ be elements of $R$ such that $Ky
\subset Kx$. Assume further that $x$ is a non zerodivisor. Let $m
> {\rm dim}_k ((K\colon {\mathfrak m}) \cap (x \colon y)/K)$ and
let $u_1, \ldots, u_m$ be units in $R$ that are not all congruent
modulo ${\mathfrak m}$. Then
\[
(K \colon {\mathfrak m})x \cap (K \colon {\mathfrak m}) (x+u_1 y)
\cap \ldots \cap (K \colon {\mathfrak m}) (x+u_m y) \subset (x
\colon y)x.
\]
\end{Lemma}
\begin{proof}
Let $\alpha$ be an element of the intersection. Then $\alpha
= s x$ with $s \in K \colon {\mathfrak m}$, but also
\[
\alpha = s_1(x+u_1 y), \ \ldots \ , \alpha = s_m(x+u_m y)
\]
with $s_i \in K \colon {\mathfrak m}$. We may assume that $s_i
\not\in K$ for every $i$, because $Ky \subset Kx$ and $K \subset x
\colon y$. On the other hand, since $\alpha \in (x)$ and $u_i$ are
units, we have $s_i \in x \colon y$. Using this we wish to show
that $s \in x \colon y$ as well.

Let `${}^{\tratto}$' denote images in $\overline{R}=R/K$. Now
$\overline{s}_1, \ldots, \overline{s}_m$, or equivalently,
$\overline{u}_1\overline{s}_1, \ldots,$
$\overline{u}_m\overline{s}_m$ are $m$ nonzero elements of the
$k$-vector space $\overline{(K \colon {\mathfrak m}) \cap (x
\colon y)}$. They are linearly dependent, and shrinking $m$ if
needed we may assume that $m$ is minimal with this property.
Obviously $m \geq 2$. Now there exist units $\lambda_1, \ldots,
\lambda_m$ in $R$ so that $\displaystyle\sum_{i=1}^m
\overline{\lambda}_i \overline{u}_i \overline{s}_i =
\overline{0}$. Notice that $\displaystyle\sum_{i=1}^m
\overline{\lambda}_i \overline{s}_i \not= \overline{0}$ because
$u_1, \ldots, u_m$ are incongruent modulo ${\mathfrak m}$. Since
$\displaystyle\sum_{i=1}^m \lambda_i u_i s_i \in K$, the element
$\displaystyle\sum_{i=1}^m \lambda_i u_i s_i y$ is in $Ky \subset
Kx$, and can be written as $\xi x$ for some $\xi \in K$. Set
$\lambda = \displaystyle \sum_{i=1}^m \lambda_i$ and multiply both
sides by $\alpha$. Rewriting $\alpha$ by means of the above
equations and cancelling $x$ we obtain
\[
\lambda s = \sum_{i=1}^m \lambda_i s_i + \xi \in (x \colon y) + K
= x \colon y.
\]
If $\lambda \in {\mathfrak m}$ then $\lambda s \in K$ since $s \in
K \colon {\mathfrak m}$, and we conclude that $\overline{0} =
\displaystyle \sum_{i=1}^m \overline{\lambda}_i \overline{s}_i$,
which is impossible. Thus $\lambda$ is a unit, and the desired
inclusion $s \in x \colon y$ follows.
\end{proof}

\smallskip

\begin{Lemma}~\label{corollary-HS}
Let $(R, {\mathfrak m})$ be a one-dimensional local
Cohen--Macaulay ring with infinite residue field. Let $I$ be an
${\mathfrak m}$-primary ideal. Assume that $J \colon I$ does not
depend on the minimal reduction $J$ of $I$. Then for every minimal
reduction $J$ of $I$
\[
{\rm core}(I) = (J \colon I) J.
\]
\end{Lemma}
\begin{proof}
By Proposition~\ref{thm1} it suffices to show the inclusion ${\rm
core}(I) \subset (J \colon I)J$. Write $J=(x)$ and $K = x \colon
I$. By assumption, for any other minimal reduction $(y)$ of $I$
one also has $K= y \colon I$. Moreover, since both $x$ and $y$ are
regular elements one can easily check that $(y \colon I)x = xy
\colon I = (x \colon I)y$. Thus $Kx=Ky$. Write $I=(y_1, \ldots,
y_n)$, where $(y_i)$ are minimal reductions of $I$. By
Lemma~\ref{lemma-HS} there exists $m \gg 0$ and minimal reductions
$J_{ij}$ of $I$, $1 \leq i \leq n$ and $0 \leq j \leq m$, so that
\[
\bigcap_{j=0}^m (K \colon {\mathfrak m}) J_{ij} \subset (x \colon
y_i)x.
\]
Since $J_{ij}=(z_{ij})$ with $z_{ij}$ $R$-regular and $Kx =
Kz_{ij}$, we have $(Kx \colon {\mathfrak m}) \cap J_{ij} \subset
(K \colon {\mathfrak m}) J_{ij}$. Thus
\begin{eqnarray*}
(Kx \colon {\mathfrak m}) \cap \bigcap_{\split{1 \leq i \leq n}{0
\leq j \leq m}} J_{ij} & \subset & \bigcap_{i=1}^n \bigcap_{j=0}^m
(K \colon {\mathfrak m})J_{ij} \subset \bigcap_{i=1}^n (x \colon y_i)x \\
& = & \left( \bigcap_{i=1}^n (x \colon y_i) \right) x = (x \colon
I) x = Kx.
\end{eqnarray*}
As the ideal $Kx$ is ${\mathfrak m}$-primary we deduce that
$\displaystyle\bigcap_{\split{1 \leq i \leq n}{0 \leq j \leq m}}
J_{ij} \subset Kx$, which gives ${\rm core}(I) \subset Kx = (J
\colon I)J$.
\end{proof}

Let $R$ be a Noetherian local ring with infinite residue field and
let $I$ be an $R$-ideal. We write ${\mathcal M}(I)$ for the set of
all minimal reductions of $I$, and define
\[
\gamma(I) = \inf \{ t \, | \, {\rm core}(I) = \bigcap_{i=1}^t J_i
\ {\rm with} \ J_i \in {\mathcal M}(I) \}.
\]
According to \cite[3.1]{CPU}, the number $\gamma(I)$ is finite for
a broad class of ideals. In particular it is finite for almost all
the ideals considered in the present paper \cite[3.2]{CPU}. Recall
that $I$ is said to be {\it equimultiple} if $\ell(I)={\rm ht}\,
I$.

\begin{Remark}\label{gamma}
Let $R$ be a local Cohen--Macaulay ring with infinite residue
field and let $I$ be an equimultiple $R$-ideal. If $R_{\mathfrak
p}$ is Gorenstein for every ${\mathfrak p} \in {\rm Min}(I)$, then
\[
\gamma(I) = \max \{ {\rm type}((R/{\rm core}(I))_{\mathfrak p}) \,
| \, {\mathfrak p} \in {\rm Min}(I) \}.
\]
\end{Remark}
\begin{proof}
One has $\gamma(I) = \max (\{ \gamma(I_{\mathfrak p}) \, |
\, {\mathfrak p} \in {\rm Min}({\rm Fitt}_g(I)) \} \cup \{ 1 \}) =
\max \{ \gamma(I_{\mathfrak p}) \, | \, {\mathfrak p} \in {\rm
Min}(I) \}$, where the first equality follows from \cite[4.9]{CPU}
and the second one is a consequence of \cite{CN}. On the other
hand if ${\mathfrak p} \in {\rm Min}(I)$, then
$\gamma(I_{\mathfrak p}) = {\rm type}(R_{\mathfrak p}/{\rm
core}(I_{\mathfrak p}))$ since every minimal reduction of
$I_{\mathfrak p}$ is an irreducible zero-dimensional $R_{\mathfrak
p}$-ideal. Finally, $R_{\mathfrak p}/{\rm core}(I_{\mathfrak p}) =
(R/{\rm core}(I))_{\mathfrak p}$ by \cite[4.8]{CPU}.
\end{proof}

\begin{Lemma}\label{free}
Let $R$ be a local Cohen--Macaulay ring and let $I$ be an
$R$-ideal. Suppose that $I$ satisfies $G_s$ and is weakly
$(s-2)$-residually $S_2$ \/ for some integer $s$. Let $K =
{\mathfrak a} \colon I$ be an $s$-residual intersection of $I$
with ${\mathfrak a} \subset I$ and $\mu({\mathfrak a}) \leq s$.
Then ${\mathfrak a}/K {\mathfrak a}$ is a free $R/K$-module of
rank $s$.
\end{Lemma}
\begin{proof} We may assume $s \geq 1$.
Let $a_1, \ldots, a_s$ be a generating sequence of
${\mathfrak a}$ and write ${\mathfrak a}_i = (a_1, \ldots,
\widehat{a}_i, \ldots, a_s)$ for $1 \leq i \leq s$. By
\cite[1.6$(${\it a}$)$]{U}, $a_1, \ldots, a_s$ can be chosen so
that ${\mathfrak a}_i \colon I$ are $(s-1)$-residual intersections
of $I$, and then by \cite[3.1]{CEU}, ${\mathfrak a}_i \colon a_i =
{\mathfrak a}_i \colon I \subset K$. Hence the coefficients of the
relations among $a_1, \ldots, a_s$ are contained in $K$.
\end{proof}

Let $(R, {\mathfrak m}, k)$ be a Noetherian local ring and let $I$
be an $R$-ideal. We say that $J_1, \ldots, J_t$ are {\it general
$s$-generated ideals in $I$} if $J_i \subset I$ are ideals with
$\mu(J_i)=s$, $J_i \otimes_R k \hookrightarrow I \otimes_R k$, and
the point $(J_1 \otimes_R k, \ldots, J_t \otimes_R k)$ lies in
some dense open subset of the product of Grassmannians
$\displaystyle\stackrel{t}{\times}{\mathbb G}(s, I \otimes_R k)$.
For an ideal $I$ of analytic spread $\ell$ we also set ${\mathcal
Q}(I) = \{ {\mathfrak p} \in V({\rm Fitt}_{\ell}(I)) \, | \, {\rm
dim}\, R_{\mathfrak p} = \ell(I_{\mathfrak p}) = \ell \}$.

\begin{Theorem}\label{converse}
Let $R$ be a local Cohen--Macaulay ring with infinite residue
field. Let $I$ be an $R$-ideal of height $g$, analytic spread
$\ell$ and minimal number of generators $n$. Suppose that $I$
satisfies $G_{\ell}$ and is universally weakly
$(\ell-1)$-residually $S_2$. Consider the following conditions:
\begin{itemize}
\item[$(${\it a}$)$]
$(J \colon I) I \subset {\rm core}(I)$ for some minimal reduction $J$ of $I$;

\item[$(${\it b}$)$]
$(J \colon I) J = (J \colon I) I = {\rm core}(I)$ for every
minimal reduction $J$ of $I$;

\item[$(${\it c}$)$]
$J \colon I$ does not depend on the minimal reduction $J$ of $I$;

\item[$(${\it d}$)$]
the reduction number of $I$ is at most $\ell-g+1$.
\end{itemize}
They relate in the following way:
\begin{itemize}
\item[$($1$)$]
Conditions $(${\it a}$)$, $(${\it b}$)$ and $(${\it c}$)$ are
equivalent, and they imply
\[
\gamma(I) \leq \ell \cdot \max \left( \{ {\rm type}((R/J \colon
I)_{\mathfrak p}) \, | \, {\mathfrak p} \in {\mathcal Q}(I) \}
\cup \{0\} \right) + 1,
\]
where $J$ is any minimal reduction of $I$.

\item[$($2$)$]
If for every ${\mathfrak p} \in {\mathcal Q}(I)$, $R_{\mathfrak
p}$ is Gorenstein and ${\rm depth} (R/I^j)_{\mathfrak p} \geq
\ell-g-j+1$ whenever $1 \leq j \leq \ell-g$, then $(${\it d}$)$
implies the other conditions, and $(${\it a}$)$-$(${\it c}$)$
imply that
\[
\gamma(I) \leq \ell \cdot {n-g \choose \ell-g+1} + 1.
\]

\item[$($3$)$]
If $\ell \geq 1$ and for every ${\mathfrak p} \in V(I)$,
$R_{\mathfrak p}$ is Gorenstein and
\[
{\rm depth} (R/I^j)_{\mathfrak p} \geq \min \{ {\rm dim}\,
R_{\mathfrak p}, \ell+1 \} -g-j+1
\]
whenever $1 \leq j \leq
\ell-g+1$, then all four conditions $(${\it a}$)$-$(${\it d}$)$
are equivalent.
\end{itemize}
\end{Theorem}
\begin{proof}
We may assume $\ell \geq 1$. Recall that by \cite[2.1$(${\it a}$)$
and 4.8]{CPU}, $({\rm core}(I))_{\mathfrak p} = {\rm
core}(I_{\mathfrak p})$ for every ${\mathfrak p} \in V(I)$.
Moreover for any minimal reduction $J$ of $I$ one has ${\rm ht}\,
J \colon I \geq \ell$ according to \cite[2.1$($\/{\it
f}\/$)$]{CPU}, and ${\rm Ass}_R(R/J \colon I) \subset {\mathcal
Q}(I)$ as well as ${\rm Ass}_R(R/J) \cap V(J \colon I) \subset
{\mathcal Q}(I)$ by \cite[3.4$(${\it a}$)$,$(${\it b}$)$]{CEU} and
\cite[2.1$(${\it e}$)$,$(${\it g}$)$]{CPU}.

We now prove $($1$)$. In light of Proposition~\ref{thm1}, to
establish the asserted equivalence it is enough to show that
$(${\it c}$)$ implies ${\rm core}(I) \subset (J \colon I)J$ for
$J$ any minimal reduction of $I$. It suffices to prove the
inclusion ${\rm core}(I) \subset (J \colon I)J$ locally at every
associated prime ${\mathfrak p}$ of the latter ideal. We may
assume that ${\mathfrak p}$ contains $J \colon I$, since the
desired inclusion is clear otherwise. Now according to
Lemma~\ref{free}, ${\rm Ass}_R(R/(J \colon I)J) \subset {\rm
Ass}_R(R/(J \colon I)) \cup {\rm Ass}_R(R/J)$. Hence by the above
we may assume that ${\mathfrak p} \in {\mathcal Q}(I)$.
Furthermore \cite[2.1$(${\it e}$)$ and 4.9]{CPU} shows that
$\gamma(I) \leq \max ( \{ \gamma(I_{\mathfrak p}) \, | \,
{\mathfrak p} \in {\mathcal Q}(I) \} \cup \{1\} )$. Finally by
\cite[2.1$(${\it a}$)$]{CPU} and Proposition~\ref{thm1}, our
assumptions are preserved as we replace $R$ by $R_{\mathfrak p}$.
Thus we may from now on assume that $R = R_{\mathfrak p}$ with
${\mathfrak p}$ in ${\mathcal Q}(I)$.

Let $a_1, \ldots, a_{\ell}$ be a generating sequence of $J$ and
write ${\mathfrak a}_i = (a_1, \ldots, \widehat{a}_i, \ldots,
a_{\ell})$, $K_i= {\mathfrak a}_i \colon I$ for $1 \leq i \leq
\ell$. By \cite[1.6$(${\it a}$)$]{U}, $a_1, \ldots, a_{\ell}$ can
be chosen so that $K_i$ are geometric $(\ell-1)$-residual
intersections.

For a fixed $i$ let `${}^{\tratto}$' denote images modulo $K_i$.
According to \cite[3.4, 3.1 and 2.4$(${\it b}$)$]{CEU},
$\overline{R}$ is a one-dimensional local Cohen--Macaulay ring,
$\overline{I}$ is a height one ideal, $K_i \cap I = {\mathfrak
a}_i$ and $\overline{J} :_{\overline{R}} \overline{I} =
\overline{J :_R I}$. We now claim that $\overline{{\rm core}(I)}
\subset {\rm core}(\overline{I})$. By \cite[4.5]{CPU}, for $t =
\gamma(\overline{I})$ general one-generated ideals $(\alpha_1),
\ldots, (\alpha_t)$ in $I$ one has ${\rm core}(\overline{I}) =
(\overline{\alpha}_1) \cap \ldots \cap (\overline{\alpha}_t)$. As
$({\mathfrak a}_i, \alpha_j)$ are reductions of $I$, we obtain
${\rm core}(I) \subset \displaystyle\bigcap_{j=1}^t ({\mathfrak
a}_i, \alpha_j)$. Therefore
\[
\overline{{\rm core}(I)} \subset \overline{\bigcap_{j=1}^t
({\mathfrak a}_i, \alpha_j)} \subset \bigcap_{j=1}^t
(\overline{\alpha}_j) = {\rm core}(\overline{I}).
\]
By Proposition~\ref{thm1}, $(J \colon I)I \subset {\rm core}(I)$.
Therefore $(\overline{J} :_{\overline{R}}
\overline{I})\overline{I} = \overline{(J :_R I)I} \subset
\overline{{\rm core}(I)} \subset {\rm core}(\overline{I})$. Thus,
again using Proposition~\ref{thm1} we deduce that $H \colon
\overline{I}$ does not depend on the minimal reduction $H$ of
$\overline{I}$. Now Lemma~\ref{corollary-HS} implies that ${\rm
core}(\overline{I}) = (\overline{J} \colon
\overline{I})\overline{J}$. Thus $\overline{{\rm core}(I)} \subset
\overline{(J \colon I)J}$, which yields ${\rm core}(I) \subset (J
\colon I)J + K_i$.

It follows that
\[
{\rm core}(I) \subset \bigcap_{1 \leq i \leq \ell} \left(\left[(J
\colon I)J +K_i \right] \cap I \right) = \bigcap_{1 \leq i \leq
\ell} \left((J \colon I)J + {\mathfrak a}_i\right).
\]
Now, using Lemma~\ref{free} we conclude that
$\displaystyle\bigcap_{1 \leq i \leq \ell} \left((J \colon I)J +
{\mathfrak a}_i\right) = (J \colon I)J$. Therefore ${\rm core}(I)
\subset (J \colon I)J$, which proves the asserted equivalence.

To establish the bound for $t = \gamma(I)$, write ${\rm core}(I) =
J_1 \cap \ldots \cap J_t$ where $J_1, \ldots, J_t$ are minimal
reductions of $I$. By $(${\it b}$)$, ${\rm core}(I) =(J_1 \colon
I) J_1$. Since the module $M = J_1/{\rm core}(I) = J_1 /(J_1
\colon I)J_1$ is Artinian, $t-1 \leq {\rm type}(M)$. On the other
hand Lemma~\ref{free} shows that ${\rm type}(M) = \ell \cdot {\rm
type}(R/J_1 \colon I)$. Finally by $(${\it c}$)$, $J_1 \colon I =
J \colon I$, completing the proof of the inequality $\gamma(I)
\leq \ell \cdot {\rm type}(R/J \colon I) +1$.

\medskip

We now prove $($2$)$. In showing that $(${\it d}$)$ implies
$(${\it a}$)$ it suffices to verify $(${\it a}$)$ locally at every
associated prime ${\mathfrak p}$ of ${\rm core}(I)$. Since by
\cite[3.1, 3.2]{CPU} ${\rm core}(I)$ is a finite intersection of
minimal reductions of $I$, we may assume that ${\mathfrak p} \in
{\mathcal Q}(I)$. If $\ell=g$ then obviously ${\rm ht}\,
I_{\mathfrak p} = g$, whereas for $\ell > g$, $\ell-{\rm ht}\,
I_{\mathfrak p} = {\rm dim}(R/I)_{\mathfrak p} \geq {\rm depth}
(R/I)_{\mathfrak p} \geq \ell-g$. Thus in either case ${\rm ht}\,
I_{\mathfrak p} = g$. Also notice that the reduction number of $I$
cannot increase upon localizing at ${\mathfrak p}$ since
${\mathfrak p} \in {\mathcal Q}(I)$. Now, replacing $R$ by
$R_{\mathfrak p}$ we may suppose that $R$ is Gorenstein of
dimension $\ell$ and ${\rm depth}\, R/I^j \geq \ell-g -j+1$ for $1
\leq j \leq \ell-g$, hence for $1 \leq j \leq \ell-g+1$. In this
situation, \cite[2.6]{U2} shows that $(${\it d}$)$ implies $(${\it
c}$)$, which in turn yields $(${\it a}$)$ by
Proposition~\ref{thm1}. Finally, the asserted inequality for
$\gamma(I)$ follows from $($1$)$ and the fact that ${\rm
type}((R/J \colon I)_{\mathfrak p}) \leq \displaystyle {n-g
\choose \ell-g+1}$ for ${\mathfrak p} \in {\mathcal Q}(I)$ by
\cite[2.9$(${\it b}$)$]{U}.

To prove $($3$)$ we need to show that $(${\it c}$)$ implies
$(${\it d}$)$. But this follows from \cite[3.3]{PU} $($see also
\cite[4.8]{U2}$)$.
\end{proof}

\begin{Remark}
In the setting of Theorem~\ref{converse}, ${\rm core}(I)$ is the
intersection of $\gamma(I)$ general $\ell$-generated ideals in $I$
which are reductions of $I$ $($see \cite[4.5]{CPU}$)$.
\end{Remark}

\begin{Remark}
Theorem~\ref{converse} is essentially a cancellation theorem
$($see \cite{CP3} for details$)$. Indeed, the fact that $(${\it
b}$)$ implies $(${\it c}$)$ as asserted in Theorem~\ref{converse}
can also be restated by saying that $(J_1 \colon I)I = (J_2 \colon
I)I$ implies $J_1 \colon I = J_2 \colon I$ for any minimal
reductions $J_1$ and $J_2$ of $I$.
\end{Remark}

\begin{Remark}\label{2.9}
Given two ideals $I$ and $J$ of $R$, the {\it coefficient ideal}
${\mathfrak a}(I,J)$ of $I$ with respect to $J$ is the largest
ideal ${\mathfrak a}$ of $R$ such that ${\mathfrak a}J={\mathfrak
a}I$. This ideal has been introduced by I.M. Aberbach and C.
Huneke \cite[2.1]{AH2} and has played an important role in
theorems of Brian\c{c}on-Skoda type with coefficients
\cite[2.7]{AH2}. Clearly, ${\mathfrak a}(I, J) \subset J \colon
I$. If, in addition, $I$ satisfies one of the equivalent
conditions $(${\it a}$)$-$(${\it c}$)$ in Theorem~\ref{converse}
then ${\mathfrak a}(I, J) = J:I$ for any minimal reduction $J$ of
$I$; in particular ${\mathfrak a}(I,J)$ does not depend on the
minimal reduction $J$ of $I$. This was already known in the case
$\ell=g$ \cite[2.1]{CP2}.
\end{Remark}

\begin{Remark}
Closely related is also the notion of the {\it adjoint} of $I$,
denoted ${\rm adj}(I)$, introduced by J. Lipman \cite{L}. If $R$
is a regular local ring essentially of finite type over a field of
characteristic zero and $I$ is an $R$-ideal with analytic spread
$\ell$, Lipman shows that ${\rm adj}(I^{\ell-1}) \subset
{\mathfrak a}(I,J)$ for any minimal reduction $J$ of $I$
\cite[2.3]{L}. Furthermore if $(R, {\mathfrak m})$ is a regular
local ring of dimension two with infinite residue field and $I$ is
an integrally closed ${\mathfrak m}$-primary ideal, he proves that
this containment is an equality \cite[3.3]{L}, and Huneke-Swanson
deduce that ${\rm core}(I) = {\rm adj}(I)\, I$ \cite[3.14]{HS}.
Combining \cite[3.5]{Hy}, Remark~\ref{2.9} and
Theorem~\ref{converse}$($2$)$, one obtains the following more
general result: Let $R$ be a regular local ring essentially of
finite type over a field of characteristic zero. Let $I$ be an
$R$-ideal of height $g$ and analytic spread $\ell$. Suppose that
$I$ satisfies $G_{\ell}$, the reduction number of $I$ is at most
$\ell-g+1$, ${\rm depth}\, R/I^j \geq {\rm dim}\, R/I - j + 1$
whenever $1 \leq j \leq \ell-g+1$, and ${\rm Proj}({\mathcal
R}(I))$ has rational singularities. Then ${\rm core}(I) = {\rm
adj}(I^{g-1})\, I$.
\end{Remark}

Huneke and Swanson were in part motivated to study the core
because the core of an integrally closed ideal in a
two-dimensional regular local ring is still integrally closed
$($see \cite[3.12]{HS}$)$. The following theorem fully generalizes
the one of Huneke and Swanson, as it is a celebrated result of O.
Zariski that integrally closed ideals in a two-dimensional regular
local ring are normal \cite[p. 385]{ZS}. We recall that $I$ is
said to be {\it normal} if $\overline{I^j}=I^j$ for all $j \geq
0$.

\begin{Theorem}\label{normality}
Let $R$ be a local Cohen--Macaulay ring with infinite residue
field and let $I$ be an $R$-ideal with analytic spread $\ell$.
Suppose that ${\rm ht}\, I > 0$, $I$ satisfies $G_{\ell}$ and
$AN_{\ell-1}^{-}$, and $I$ is universally weakly
$(\ell-1)$-residually $S_2$. Further assume that one of the
equivalent conditions $(${\it a}$)$-$(${\it c}$)$ of\/ {\rm
Theorem~\ref{converse}} holds. Let $J$ be a minimal reduction of
$I$.
\begin{itemize}
\item[$($1$)$]
The ${\mathcal R}(I)$-module $(J \colon I){\mathcal R}(I)$
satisfies $S_2$. This module is a maximal Cohen--Macaulay
${\mathcal R}(I)$-module if $I$ satisfies $AN_{\ell}$.

\item[$($2$)$]
Assume in addition that $R$ is a normal ring and $I$ is a normal
ideal. Then the ideal $(J \colon I)I^j$ is integrally closed for
every $j \geq 0$. In particular $J \colon I$ and ${\rm core}(I)$
are integrally closed.
\end{itemize}
\end{Theorem}
\begin{proof} Write $d = {\rm dim}\, R$ and $K = J \colon I$.

To prove $($1$)$ it suffices to show that as an ${\mathcal
R}(J)$-module, $K {\mathcal R}(I)$ satisfies $S_2$ or is maximal
Cohen--Macaulay, respectively. Condition $(${\it b}$)$ of
Theorem~\ref{converse} implies that $K {\mathcal R}(I) = K
{\mathcal R}(J)$. According to \cite[2.1$($\/{\it f}\/$)$]{CPU},
${\rm ht}\, K \geq \ell$. Therefore $J$ satisfies $G_{\infty}$ and
$AN_{\ell-1}^{-}$ by \cite[1.12]{U}, and then $J$ has the `sliding
depth property' by \cite[1.8$(${\it c}$)$]{U}. Now \cite[9.1]{HSV}
shows that ${\mathcal R}(J) \simeq {\rm Sym}_R(J)$ is a
Cohen--Macaulay ring, necessarily of dimension $d+1$. Thus we may
assume that $K \not= R$. Now $K$ is unmixed of height $\ell$
according to \cite[1.7$(${\it a}$)$]{U}, and if $I$ satisfies
$AN_{\ell}$ then $R/K$ is Cohen--Macaulay. On the other hand
Lemma~\ref{free} gives that $J/KJ$ is a free $R/K$-module of rank
$\ell$, and therefore ${\mathcal R}(J)/K {\mathcal R}(J) \simeq
{\rm Sym}_{R/K}(J/KJ)$ is a polynomial ring over $R/K$ in $\ell$
variables. Thus ${\mathcal R}(J)/K {\mathcal R}(J)$ is $S_1$ and
equidimensional of dimension $d$, and if $I$ satisfies $AN_{\ell}$
then this ring is Cohen--Macaulay. As ${\mathcal R}(J)$ is
Cohen--Macaulay of dimension $d+1$ we deduce that the ${\mathcal
R}(J)$-module $K{\mathcal R}(J)$ satisfies $S_2$ or is maximal
Cohen--Macaulay, respectively.

We now prove $($2$)$. By part $($1$)$, $K {\mathcal R}(I)$ is a
divisorial ideal of the normal domain ${\mathcal R}(I)$. Such an
ideal is necessarily integrally closed. Since $K {\mathcal R}(I) =
\displaystyle\bigoplus_{j \geq 0} KI^j t^j$ we conclude that the
$R$-ideals $KI^j$ are integrally closed as well $($for a similar
argument see \cite[the proof of 1.5]{HH}$)$. \end{proof}

\section{Examples}

In this section, we illustrate Theorems~\ref{converse}
and~\ref{normality} with a sequence of examples.

\begin{Example}
{\rm {\it Good} ideals in the sense of S. Goto, S.-I. Iai and
K.-I. Watanabe \cite{GIW} -- by which we mean ${\mathfrak
m}$-primary ideals $I$ of a local Gorenstein ring $(R, {\mathfrak
m})$ such that $I = J \colon I$ and $I^2=JI$ for some complete
intersection ideal $J$ -- are examples of ideals for which ${\rm
core}(I) = I^2$. }
\end{Example}

\begin{Example}\label{2rlr}
{\rm Let $R$ be a two-dimensional regular local ring with infinite
residue field. It is well known that all powers of the maximal
ideal ${\mathfrak m}=(x,y)$ have reduction number at most one.
From Theorem~\ref{converse}$($2$)$ we obtain
\[
{\rm core}({\mathfrak m}^j) = ((x^j, y^j) \colon {\mathfrak m}^j)
{\mathfrak m}^j = {\mathfrak m}^{j-1}{\mathfrak m}^j = {\mathfrak
m}^{2j-1}.
\]
In particular, $\gamma({\mathfrak m}^j) = 2j-1$ by
Remark~\ref{gamma}. This shows that the bound in
Theorem~\ref{converse}$($2$)$ is sharp and that there is no upper
bound for $\gamma(I)$ as $I$ varies among the ideals of $R$. If
the restriction on the dimension is dropped, one can still compute
the core of ${\mathfrak m}^j$, although Theorem~\ref{converse} may
no longer apply, see Proposition~\ref{proposition4.4}. }
\end{Example}

\begin{Example}{\rm
Let $R = k[x, y, z]_{(x, y, z)}$ with $k$ a field of
characteristic zero and $x, y, z$ variables, and let ${\mathfrak
m}$ denote the maximal ideal of $R$. Consider the Gorenstein
$R$-ideal of height three
\[
I = (x^2-y^2+xz, xy+xz-yz, xz-2yz+z^2, y^2+yz-z^2, z^2-2yz),
\]
and let $J$ be the $R$-ideal generated by the first three
generators of $I$. The ideal $I$ has reduction number $2$. Hence
it fails to satisfy conditions $(${\it a}$)$-$(${\it d}$)$ of
Theorem~\ref{converse}. However, by \cite[4.5]{CPU} the core of
$I$ is the intersection of homogeneous minimal reductions of $I$.
So it will definitely contain ${\mathfrak m}^4$. On the other
hand, taking $10$ general homogeneous minimal reductions $J_1,
\ldots, J_{10}$ of $I$ a calculation using the computer algebra
system {\tt Macaulay} shows that $J_1 \cap \ldots \cap J_{10} =
{\mathfrak m}^4$. Therefore ${\rm core}(I)={\mathfrak m}^4$. }
\end{Example}

In the case of perfect ideals of height two or perfect Gorenstein
ideals of height three one can compute the core explicitly from a
matrix presenting $I$.

\begin{Corollary}\label{explicit1}
Let $R$ be a local Gorenstein ring with infinite residue field.
Let $I$ be a perfect $R$-ideal of height two and analytic spread
$\ell$, satisfying $G_{\ell}$. Let $\varphi$ be a matrix with $n$
rows presenting $I$. Then the following conditions are equivalent:
\begin{itemize}
\item[$(${\it a}$)$]
the reduction number of $I$ is at most $\ell-1$;

\item[$(${\it b}$)$]
${\rm core}(I) = I_{n-\ell}(\varphi) \cdot I$.
\end{itemize}
If any of these conditions holds and if $R$ and $I$ are normal,
then $I_{n-\ell}(\varphi)$ and ${\rm core}(I)$ are integrally
closed.
\end{Corollary}
\begin{proof} The result follows from \cite{BE} and
Theorems~\ref{converse}$($3$)$ and~\ref{normality}$($2$)$.
\end{proof}

\begin{Corollary}\label{corHS} $(${\rm see \cite[3.9]{HS}}$)$
Let $R$ be a two-dimensional regular local ring with infinite
residue field. Let $I$ be an integrally closed $R$-ideal and let
$\varphi$ be a matrix with $n$ rows presenting $I$. Then
\[
{\rm core}(I) =I_{n-2}(\varphi) \cdot I.
\]
Furthermore $I_{n-2}(\varphi)$ and ${\rm core}(I)$ are integrally
closed.
\end{Corollary}
\begin{proof} We may assume that $I$ is primary to the maximal ideal of
$R$. The ideal $I$ has reduction number at most one by
\cite[5.5]{LT} and is normal according to \cite[p. 385]{ZS}. Now
the assertion follows from Corollary~\ref{explicit1}.
\end{proof}

\begin{Corollary}\label{explicit2}
Let $R$ be a local Gorenstein ring with infinite residue field.
Let $I$ be a perfect Gorenstein $R$-ideal of height three,
analytic spread $\ell$ and minimal number of generators $n$,
satisfying $G_{\ell}$. Let $\varphi$ be a matrix with $n$ rows
presenting $I$. Then the following conditions are equivalent:
\begin{itemize}
\item[$(${\it a}$)$]
the reduction number of $I$ is $\ell-2$;

\item[$(${\it b}$)$]
${\rm core}(I) = I_1(\varphi) \cdot I$.
\end{itemize}
If any of these conditions holds and if $R$ and $I$ are normal,
then  $I_1(\varphi)$ and ${\rm core}(I)$ are integrally closed.
\end{Corollary}
\begin{proof}
The result follows from \cite[5.6]{U2}, \cite[5.5]{JU} and
Theorems~\ref{converse}$($3$)$ and~\ref{normality}$($2$)$.
\end{proof}

\begin{Corollary}\label{explicit3}
Let $R$ be a local Gorenstein ring with infinite residue field.
Let $I$ be an $R$-ideal of height $g$, analytic spread $\ell$,
minimal number of generators $n = \ell+1 \geq 2$, satisfying
$G_{\ell}$ and ${\rm depth}\, R/I^j \geq {\rm dim}\, R/I -j+1$ for
$1 \leq j \leq \ell-g+1$. Let $\varphi$ be a matrix with $n$ rows
presenting $I$. Then  the following conditions are equivalent:
\begin{itemize}
\item[$(${\it a}$)$]
the reduction number of $I$ is $\ell-g+1$;

\item[$(${\it b}$)$]
${\rm core}(I) = I_1(\varphi) \cdot I$.
\end{itemize}
If any of these conditions holds and if $R$ and $I$ are normal,
then  $I_1(\varphi)$ and ${\rm core}(I)$ are integrally closed.
\end{Corollary}
\begin{proof}
The result follows from \cite[5.5]{JU} and
Theorems~\ref{converse}$($3$)$ and~\ref{normality}$($2$)$.
\end{proof}

\begin{Remark}
In the setting of Corollaries~\ref{explicit1},~\ref{explicit2}
and~\ref{explicit3}, \cite[2.6]{U2} $($or
Theorem~\ref{converse}$($2$))$ provides additional simplification
in the computation of the core of $I$. In the case of
Corollary~\ref{explicit1}$(${\it b}$)$ we have that after
elementary row operations, $I_{n-\ell}(\varphi)$ is generated by
the $n-\ell$ by $n-\ell$ minors of the matrix consisting of the
last $n-\ell$ rows of $\varphi$. In the case of
Corollaries~\ref{explicit2}$(${\it b}$)$
and~\ref{explicit3}$(${\it b}$)$ we have that after elementary row
operations, $I_1(\varphi)$ is generated by the entries of the last
row of $\varphi$.
\end{Remark}

The next example is an illustration of Corollary~\ref{explicit3}.
In this situation the integral closedness of the ideal alone
suffices to prove the one of the core.

\begin{Example}\label{linear}
Let $R = k[x_1, \ldots, x_d]_{(x_1, \ldots, x_d)}$ with $k$ an
infinite field and $x_1, \ldots, x_d$ variables, and let
${\mathfrak m}$ denote the maximal ideal of $R$. Let $I$ be an
$R$-ideal of height $g$ and minimal number of generators $d+1 \geq
2$, satisfying $G_d$ and ${\rm depth}\, R/I^j \geq {\rm dim}\,
R/I-j+1$ for $1 \leq j \leq d-g$. Assume that $I$ has a
presentation matrix whose entries are linear forms in $k[x_1,
\ldots, x_d]$. Then ${\rm core}(I) = {\mathfrak m}I$. If $I$ is
integrally closed then so is ${\rm core}(I)$.
\end{Example}
\begin{proof}
Let $\varphi$ be a matrix of linear forms presenting $I$. Notice
that after elementary row operations, the entries of the last row
of $\varphi$ generate ${\mathfrak m}$. Now by \cite[5.1]{U2},
$\ell(I)=d$ and condition $(${\it a}$)$ of
Corollary~\ref{explicit3} holds. Since $I_1(\varphi) = {\mathfrak
m}$, the corollary shows that ${\rm core}(I) = {\mathfrak m}I$.
The ideal $I$ is generated by forms in $k[x_1, \ldots, x_d]$ of
the same degree, say $s$. Therefore ${\mathfrak m}I = {\mathfrak
m}^{s+1} \cap I$, and it follows that ${\rm core}(I)$ is
integrally closed if $I$ is.
\end{proof}

Given the previous example, one may hope that a more general
analogue of Theorem~\ref{normality}$($2$)$ -- where the normality
of the ideal $I$ is replaced by $I$ being integrally closed --
holds. The following example however shows that this is not true
in general.

\begin{Example}\label{bad-pfaff}
Let $R = k[x, y, z, w]_{(x,y,z,w)}$ with $k$ a field of
characteristic zero. Let $I$ be the $R$-ideal generated by the
four by four Pfaffians of the five by five alternating matrix
\[
\varphi=\left(
\begin{array}{ccccc}
0   & -x^2 & -y^2 & -z^2 & -w^2 \\
x^2 & 0    & -w^2 & y^2  & -z^2 \\
y^2 & w^2  & 0    & -x^2 & -x^2 \\
z^2 & -y^2 & x^2  & 0    & -y^2 \\
w^2 & z^2  & x^2  & y^2  & 0
\end{array}
\right).
\]
The ideal $I$ is perfect Gorenstein of height three. Since $I$ is
generically a complete intersection and $I_1(\varphi)$ is
generated by the entries of the last row of $\varphi$, by
\cite[5.1]{U2} $I$ satisfies condition $(${\it a}$)$ of
Corollary~\ref{explicit2}. Thus one has that ${\rm core}(I) =
(x^2, y^2, z^2, w^2) \cdot I$. A computation using the computer
algebra system {\tt Macaulay} shows that $I = \surd I$. Hence the
ideal $I$ is integrally closed. However, ${\rm core}(I)$ is not an
integrally closed ideal: For example, the element
$xw(x^4+y^4+z^2w^2) \not\in {\rm core}(I)$, but it is integral
over ${\rm core}(I)$.
\end{Example}

\section{Symbolic powers of ideals}

Let $I$ be an ideal in a Noetherian ring $R$, let $W$ be the
complement in $R$ of the union of all the associated primes of
$I$, and let $j \geq 1$ be an integer. We recall that the $j^{\rm
th}$ {\it symbolic power} $I^{(j)}$ of $I$ is the preimage of $I^j
R_W$ in $R$. In this section we study the core of links of
symbolic powers of ideals. We also compute the core of powers of
complete intersections. Finally, we apply the latter result to the
case of even symbolic powers of self-linked ideals.

\begin{Proposition}\label{symbolic}
Let $R$ be a local Gorenstein ring with infinite residue field.
Let $H$ be an unmixed $R$-ideal of height $g \geq 2$ which is
generically a complete intersection. Let $\underline{\mathbf
\alpha}= \alpha_1, \ldots , \alpha_g$ be a regular sequence
contained in $H^{(j)}$ and set $I = (\underline{\mathbf \alpha})
\colon H^{(j)} = (\underline{\mathbf \alpha}) \colon H^j$ for some
$j \geq 1$. If $g=2$ or $j=1$, assume that at least two of the
$\alpha_i$'s are contained in $H^{(j+1)}$. Then
\[
{\rm core}(I) = H^{(j)} I.
\]
\end{Proposition}
\begin{proof}
From \cite[4.1]{PU} we have $I^2 = (\underline{\mathbf
\alpha}) I$. Hence Theorem~\ref{converse}$($2$)$ shows that ${\rm
core}(I) = ((\underline{\mathbf \alpha}) \colon I) I = H^{(j)} I$.
\end{proof}

Sometimes the core of $($symbolic$)$ powers can be computed
directly even if the reduction number is not the expected one:

\begin{Proposition}\label{proposition4.4}
Let $R$ be a Noetherian local ring with infinite residue field.
Let $I$ be an $R$-ideal of height $g \geq 1$ generated by a
regular sequence. Then for any $j \geq 1$,
\[
{\rm core}(I^j) = I^{gj-g+1}.
\]
If $R$ is Gorenstein one has $\gamma(I^j) = \displaystyle{gj-1
\choose g-1}$.
\end{Proposition}
\begin{proof}
The second assertion is an immediate consequence of
Remark~\ref{gamma}. To prove the first one write ${\mathfrak m}$
for the maximal ideal of $R$, $k = R/{\mathfrak m}$, $A = R/I$,
${\mathfrak n} = {\mathfrak m}/I$, ${\mathcal G} = {\rm gr}_I(R)$,
and ${\mathcal F} = {\rm gr}_I(R) \otimes_R k$. Notice that
${\mathcal G} \simeq A[x_1, \ldots, x_g]$ and ${\mathcal F} \simeq
k[x_1, \ldots, x_g]$ are polynomial rings in $g$ variables with
irrelevant ideals ${\mathcal G}_{+} = (x_1, \ldots, x_g){\mathcal
G}$ and ${\mathcal F}_{+} = (x_1, \ldots, x_g){\mathcal F}$,
respectively.

Let $J$ be any minimal reduction of $I^j$ and set $J' = J +
{\mathfrak m}I^j /{\mathfrak m}I^j \subset [{\mathcal F}]_j$. The
equality $(I^j)^{r+1} = J(I^j)^r$ for $r \gg 0$ implies that
$\sqrt{J' {\mathcal F}} = {\mathcal F}_{+}$. Thus $J'{\mathcal F}$
is a complete intersection ${\mathcal F}$-ideal minimally
generated by $g$ forms of degree $j$. Computing the socle degree
of ${\mathcal F}/J'{\mathcal F}$, we see that $({\mathcal
F}_{+})^{gj-g+1} \subset J'{\mathcal F}$. Thus
$I^{gj-g+1}=JI^{gj-g-j+1}$ by Nakayama's Lemma. This proves the
inclusion $I^{gj-g+1} \subset {\rm core}(I^j)$.

To show the opposite containment, let $u \in R \setminus
I^{gj-g+1}$. We need to find a minimal reduction of $I^j$ that
does not contain $u$. For this we may suppose $gj-g \geq j$.
Furthermore, as ${\mathcal G}$ is a polynomial ring we may replace
$u$ by a suitable multiple to assume that $u \in I^{gj-g}
\setminus I^{gj-g+1}$. Write $f$ for the image of $u$ in
$[{\mathcal G}]_{gj-g}$. After a change of generators of $I$,
which amounts to a change of variables in the polynomial ring
${\mathcal G}$, we may assume that $f \in {\mathfrak n}^s{\mathcal
G}$ for some $s$, but that the coefficient of $x_1^{gj-g}$ in $f$
does not belong to ${\mathfrak n}^{s+1}$. For $\alpha \in A$
consider the ${\mathcal G}$-ideals $H_{\alpha} = (\{ x_i^j-\alpha
x_i^{j-1}x_{i+1} \, | \, 1 \leq i \leq g-1 \} \cup \{ x_g^j \})$.
Since $\sqrt{H_{\alpha}} = \sqrt{{\mathcal G}_{+}}$ and
$H_{\alpha}$ are generated by $g$ forms of degree $j$, it follows
that $H_{\alpha}$ lift to minimal reductions $J_{\alpha}$ of
$I^j$. By \cite[2.7]{VV} one has $J_{\alpha} \cap I^{gj-g} =
J_{\alpha} I^{gj-g-j}$ because $H_{\alpha}$ is generated by a
regular sequence on ${\mathcal G}$. Thus if $u \in J_{\alpha}$
then $u \in J_{\alpha}I^{gj-g-j}$, which gives $f \in H_{\alpha}$.
Hence it suffices to show that $f \not\in H_{\alpha}$ for some
$\alpha \in A$.

To do so we map $A$ onto a ring $\overline{A}$ so that
$\overline{{\mathfrak n}}^{s+1} = 0$, $\overline{{\mathfrak n}}^s
\simeq k$, and the coefficient of $x_1^{gj-g}$ in $\overline{f}$
is not zero, where `${}^{\tratto}$' denotes images in
$\overline{A}$ and $\overline{A}[x_1, \ldots, x_g]$, respectively.
Recall that $f \in \overline{{\mathfrak n}}^s[x_1, \ldots, x_g]$.
Reverting to our original notation we write $A$ instead of
$\overline{A}$. Now $f = ah$, with $a \not= 0$ in $A$ and $h$ a
homogeneous polynomial of degree $gj-g$ that is monic in $x_1$.

Set $B = A[x_1, \ldots, x_g]/H_{\alpha}$ and notice that $B$ is
flat, hence free, over $A$ because the minimal generators of
$H_{\alpha}$ form a regular sequence on the factor ring $k[x_1,
\ldots, x_g]$. Thus if $h \not\in {\mathfrak n}[x_1, \ldots, x_g]
+ H_{\alpha}$, then the image of $h$ in $B$ is linearly
independent over $A$, and hence the image of $f=ah$ in $B$ is
nonzero, as desired. Finally, to show that $h \not\in {\mathfrak
n}[x_1, \ldots, x_g]+H_{\alpha}$ for some $\alpha \in A$, we
replace $A$ by its residue field $k$. Write $M =
\displaystyle\prod_{i=1}^g x_i^{j-1}$. For every $\alpha \in k$,
$x_1^{gj-g} \equiv \alpha^{{}^{{g \choose 2}(j-1)}}M\mod
H_{\alpha}$ and any other monomial of degree $gj-g$ in $k[x_1,
\ldots, x_g]$ is either in $H_{\alpha}$ or is congruent to
$\alpha^l M$ for some $l < {g \choose 2}(j-1)$. From this we
conclude, first, that $M \not\in H_{\alpha}$ since $[B]_{gj-g}
\not= 0$. Second, $h \equiv q(\alpha) M\mod H_{\alpha}$, where
$q(Z)$ is a nonzero polynomial in $k[Z]$. Thus $q(\alpha) M
\not\in H_{\alpha}$ for some $\alpha \in k$, and then $h \not\in
H_{\alpha}$ as required.
\end{proof}

\begin{Remark}
The above proof can be simplified substantially for rings of
residue characteristic zero. In conjunction with \cite[4.5]{CPU},
Proposition~\ref{proposition4.4} says in particular that if $k$ is
an infinite field and $R = k[x_1, \ldots, x_d]$ is a polynomial
ring with irrelevant maximal ideal ${\mathfrak m}$, then
$\displaystyle{d\,j-1 \choose d-1}$ general ideals generated by
$d$ forms of degree $j$ intersect in ${\mathfrak m}^{dj-d+1}$.
\end{Remark}

The ideals of the next corollary have been extensively studied by
S. Kleiman and B. Ulrich in \cite{KU}.

\begin{Corollary}\label{even-core}
Let $R$ be a local Cohen--Macaulay ring with infinite residue
field and let $I$ be a perfect $R$-ideal of height $2$. Assume
that $I$ is generically a complete intersection and that $I =
(\Delta, \alpha) \colon I$ for some $\Delta \in I^{(2)}$ and some
$\alpha \in I$. Then for any $j \geq 1$,
\[
{\rm core}(I^{(2j)}) = I^{(4j-1)}.
\]
Furthermore, if $R$ is Gorenstein one has $\gamma(I^{(2j)})=4j-1$.
\end{Corollary}
\begin{proof}
According to \cite[3.3$($2$)$]{KU}, $I^{(2j)}$ is
equimultiple. Thus by \cite[3.1]{CPU}, ${\rm core}(I^{(2j)})$ is a
finite intersection of complete intersection ideals that are
reductions of $I^{(2j)}$. Hence it suffices to check the equality
${\rm core}(I^{(2j)}) = I^{(4j-1)}$ locally at every minimal prime
${\mathfrak p}$ of $I$. Recall that $({\rm
core}(I^{(2j)}))_{\mathfrak p} = {\rm core}((I^{(2j)})_{\mathfrak
p})$ by \cite[4.8]{CPU}. Now $I_{\mathfrak p}$ is a complete
intersection, and the powers and symbolic powers of $I_{\mathfrak
p}$ coincide. The asserted equality then follows from
Proposition~\ref{proposition4.4}. Finally, $\gamma(I^{(2j)})=4j-1$
by Remark~\ref{gamma}.
\end{proof}

\section{A conjecture about the core of ideals}

We conclude the paper by analyzing a very general formula, which
involves a fairly broad class of ideals. This class includes, for
example, all ${\mathfrak m}$-primary $($or more generally
equimultiple$)$ ideals. The thrust of Conjecture~\ref{conj} below
is to move away from the balancedness of the ideal $I$ -- which is
the main restriction in Theorem~\ref{converse} -- as much as
possible.

\begin{Conjecture}\label{conj}
Let $R$ be a local Cohen--Macaulay ring with infinite residue
field. Let $I$ be an $R$-ideal of analytic spread $\ell \geq 1$
that satisfies $G_{\ell}$ and is weakly $(\ell-1)$-residually
$S_2$. Let $J$ be a minimal reduction of $I$ and let $r$ denote
the reduction number of $I$ with respect to $J$. Then
\[
{\rm core}(I) = (J^r \colon I^r)I = (J^r \colon I^r)J = J^{r+1}
\colon I^r.
\]
\end{Conjecture}

The conjectured formula has solid theoretical foundation and an
extensive computer evidence in its support. Moreover it is sharp,
in the sense that \cite[4.11]{CPU} provides a counterexample if
the previous assumptions are relaxed. However, the situation is
now more complicated than the one encountered in
Theorem~\ref{converse} as the reduction number $r$ in
Conjecture~\ref{conj} may very well depend on the chosen minimal
reduction $J$ of $I$.

\begin{Example}
Let $R=k[x,y]_{(x,y)}$ with $k$ a field of characteristic zero.
Consider the $R$-ideals $I = (x^7, x^6y, x^2y^5, y^7)$, $J =(x^7,
y^7)$ and $H = (x^7, x^6y+y^7)$. We used the computer algebra
system {\tt Macau\-lay} and checked that $J$ and $H$ are minimal
reductions of $I$ with $r_J(I)=4$ and $r_H(I)=3$, respectively.
Using the algorithm we designed as a corollary of \cite[5.4]{CPU}
we obtained that ${\rm core}(I) = (x, y)^{13}$, and we also
verified that $(x, y)^{13} = (J^4 \colon I^4)I = (H^3 \colon
I^3)I$.
\end{Example}

The next result shows that Conjecture~\ref{conj} is also
consistent with the findings of Theorem~\ref{converse}, at least
in the case where the ring $R$ is Gorenstein and the powers of the
ideal $I$ satisfy sliding depth conditions.

\begin{Proposition}\label{agreement}
Let $R$ be a local Gorenstein ring with infinite residue field.
Let $I$ be an $R$-ideal of height $g \geq 2$, analytic spread
$\ell$, and reduction number $r$. Suppose that $I$ satisfies
$G_{\ell}$, that ${\rm depth}\, R/I^j \geq {\rm dim}\, R/I-j+1$
whenever $1 \leq j \leq \ell-g+1$, and that $r \leq \ell-g+1$.
Then for any minimal reduction $J$ of $I$,
\[
{\rm core}(I) = (J^r \colon I^r) I = (J^r \colon I^r) J = J^{r+1}
\colon I^r.
\]
\end{Proposition}
\begin{proof}
We may assume $r \geq 1$. Choose general generators $a_1, \ldots,
a_{\ell}$ for $J$ as in \cite[1.6$(${\it b}$)$]{U}. By
Theorem~\ref{converse}$($2$)$ we know that ${\rm core}(I) = (J
\colon I) I = (J \colon I) J$. According to \cite[3.4 and
4.6]{JU}, the Rees algebra ${\mathcal R}(I)$ is Cohen--Macaulay
and the reduction number of $I$ with respect to $J$ is $r$. But
then \cite[3.4]{Hy} gives the equality $J \colon I = JI^{r-1}
\colon I^r = J^r \colon I^r$. $($Observe that $JI^{r-1} \colon I^r
= J \colon I$ also follows from \cite[2.3]{CP3}; indeed, $JI^{r-1}
\colon I^r = (JI^{r-1} \colon I^{r-1}) \colon I = J \colon I$.$)$
Thus, we are left to show that
\[
(J^r \colon I^r) J = J^{r+1} \colon I^r.
\]
Clearly, $(J^r \colon I^r) J \subset J^{r+1} \colon I^r$.
Conversely, pick $\alpha \in J^{r+1} \colon I^r$. As in the proof
of Theorem~\ref{normality} one sees that the associated graded
ring ${\mathcal G}$ of $J$ is Cohen--Macaulay. Hence ${\rm
grade}\, {\mathcal G}_{+} > 0$, which yields $J^{r+1} \colon I^r
\subset J^{r+1} \colon J^r = J$. Thus we can write
\[
\alpha = \lambda_1 a_1 + \ldots + \lambda_{\ell} a_{\ell},
\]
for some $\lambda_i \in R$. To conclude, it suffices to show that
$\lambda_i \in JI^{r-1} \colon I^r = J^r \colon I^r$ for $1 \leq i
\leq \ell$ . It is enough to prove the claim for $i = \ell$. Let
$c$ be any element of $I^r$ and compare the equations
\[
\alpha c = \lambda_1 c a_1 + \ldots + \lambda_{\ell-1} c
a_{\ell-1} + \lambda_{\ell} c a_{\ell} \qquad \alpha c = \gamma_1
a_1 + \ldots + \gamma_{\ell-1} a_{\ell-1} + \gamma_{\ell}
a_{\ell},
\]
where $\gamma_i \in J^r$. The latter equation arises from the fact
that $\alpha \in J^{r+1} \colon I^r$. In conclusion we have that
\[
\lambda_{\ell} c -\gamma_{\ell} \in ((a_1, \ldots, a_{\ell-1})
\colon a_{\ell}) \cap I^r = (a_1, \ldots, a_{\ell-1})I^{r-1},
\]
by \cite[2.5$(${\it b}$)$]{JU}. Thus $\lambda_{\ell} c \in
JI^{r-1}$ for any $c \in I^r$ or, equivalently, $\lambda_{\ell}
\in JI^{r-1} \colon I^r = J^r \colon I^r$ as desired.
\end{proof}

The ideals of Proposition~\ref{proposition4.4} provide another
class of examples for which Conjecture~\ref{conj} holds.

\end{document}